%Paper: alg-geom/9403002
%From: bernd@math.cornell.edu (Bernd Sturmfels)
%Date: Tue, 1 Mar 1994 17:07:28 -0500

\magnification=\magstep1
\baselineskip=14pt
\def \Box {\hbox{}\nobreak \vrule width 1.6mm height 1.6mm
depth 0mm  \par \goodbreak \smallskip}

\vskip .4cm

\centerline{\bf INTERSECTION THEORY ON TORIC VARIETIES}

\vskip .5cm

{\baselineskip=12pt

\centerline{ {\bf William Fulton}\footnote{${}^1$}
{Supported in part by the N.S.F.}}
\centerline{ Department of Mathematics}
\centerline{University of Chicago}
\centerline{Chicago, Illinois 60637}
\centerline{\tt fulton@math.uchicago.edu}

\vskip .2cm

\centerline{and}

\vskip .2cm

\centerline{ {\bf Bernd Sturmfels}\footnote{${}^2$}{Supported in
part by the N.S.F. and the David and Lucile Packard Foundation}}
\centerline{ Department of Mathematics}
\centerline{Cornell University }
\centerline{ Ithaca, New York 14853}
\centerline{\tt bernd@math.cornell.edu}

}

\vskip 1cm

\midinsert
\narrower \narrower
{
\centerline{\sl Abstract}
\noindent The operational Chow cohomology classes of a complete toric
variety are identified with certain functions, called Minkowski weights,
on the corresponding fan.  The natural product of Chow
cohomology classes makes the Minkowski weights into a commutative
ring; the product is computed by a displacement in the lattice, which
corresponds to a deformation in the toric variety.  We show that, with
rational coefficients, this ring embeds in
McMullen's polytope algebra, and
that the polytope algebra is the direct limit of these Chow rings, over all
compactifications of a given torus.  In the nonsingular case, the Minkowski
weight corresponding to the Todd class is related to a certain Ehrhart
polynomial.}
\endinsert

\vfill
\eject

\centerline{\bf Introduction }

\vskip .1cm
\noindent
Any algebraic variety  $X$  has Chow ``homology'' groups  $A_*(X)
=  \oplus_k A_k(X)$  and Chow ``cohomology'' groups
$A^*(X) = \oplus_k A^k(X)$. The latter are the operational groups
defined in [7].  These cohomology  groups have a natural ring
structure, written with a cup product, and   $A_*(X)$  is a module
over $ A^*(X)$,  written with a cap product.  These  satisfy
functorial properties similar to homology and cohomology groups
in topology, and vector bundles have Chern classes in $A^*(X)$. If
$ X$  is complete, one has a Kronecker duality homomorphism
$$  {\cal D}_X \,\,:\,\, A^k(X) \quad \rightarrow  \quad
Hom(A_k(X),{\bf Z}), $$
that takes  $c$  to the map  $a \mapsto
deg(c \cap a)$.   For general varieties, even if nonsingular,
these groups are very  large and hard to compute, and ${\cal D}_X$
is far from being an isomorphism.

If  $X$  is a toric variety, however, $ A_k(X)$  is generated by the orbit
closures $ V(\sigma)$,  as  $\sigma $  varies over the cones of
codimension $ k$  in  the fan $\Delta$  in a lattice $ N$
corresponding to $ X$.  The relations are given
by the divisors of torus invariant rational functions on
$V(\tau)$,  for  $\tau$  a cone of  codimension  $k+1$.
This gives an explicit presentation of the
Chow groups (Proposition 1.1).  In addition, if  $X$  is complete, the
Kronecker duality homomorphism ${\cal D}_X$
is an isomorphism.  This identifies Chow
cohomology classes with certain functions on
the set of cones in $ \Delta $;
a function $c$ corresponds to a class in  $A^k(X)$
if it satisfies the linear equations in
formula (2.1) below.
We call these functions {\it Minkowski weights}.

The Chow homology of a complete toric variety can have torsion,
and  some of these groups can vanish (Example 1.3).  From the
duality isomorphism ${\cal D}_X$, we see that the Chow cohomology
groups are always torsion  free.  We include a description of
these groups for the toric varieties  corresponding to
hypersimplices (Proposition 2.6).

There is a canonical homomorphism from  $ Pic(X) $ to $ A^1(X)$,  which
is not always an isomorphism, but is an isomorphism when $ X$  is a
complete toric variety.  This is a general result of Brion [2]
for spherical  varieties, given a simple proof here in terms of
Minkowski weights  (Corollary 2.4) for the toric case.

The ring structure on $ A^*(X)$  makes the Minkowski weights into
a  commutative, associative ring.  We prove the following formula
for the  product.

\proclaim Theorem.  For  $c  \in  A^p(X)$, $ {\tilde c}  \in
A^q(X)$,  the product  $ c \cup
{\tilde c}  $ in $ A^{p+q}(X)$
  is given by the Minkowski weight that assigns to a
cone  $\gamma$  of codimension  $p+q $ the value
$$   ( c \cup {\tilde c})(\gamma)  \quad = \quad  \sum
m_{\sigma,\tau}^\gamma  \cdot  c(\sigma) \cdot {\tilde c}(\tau). $$
The sum is over a certain set of cones  $\sigma $  and  $\tau$
of codimension  $p $ and $ q$  that contain  $\gamma$,
determined by the choice of a generic vector
$v$  in $ N$:  $\sigma$  and $ \tau$  appear when
$\sigma + v$  meets  $\tau$.  The
coefficient $m_{\sigma,\tau}^{\gamma}$  is the index
$ [N : N_\sigma + N_\tau]$
where $N_\sigma := {\bf Z}(N \cap \sigma)$ and
$N_\tau := {\bf Z}(N \cap \tau)$.

\eject

The present paper is a sequel to [8], where the general results
stated before the theorem were proved for varieties on which
a solvable linear  algebraic group acts with a finite number of
orbits.  In particular, it was shown in [8] that the cup product is
given by a formula as in the theorem, where the  coefficients  $\,
m_{\sigma,\tau}^\gamma$  are obtained by expressing the diagonal
class of  $V(\gamma)$  as a linear combination of the classes of
$ V(\sigma ) \times V(\tau)$. What is new here is a
much more explicit combinatorial description of the numbers
$m_{\sigma,\tau}^\gamma$ for toric varieties.
In order to prove the theorem, we carry out an
explicit  rational deformation, flowing in the direction
corresponding to the vector   $v$ (Section 3).

The Chow rings $A^*(X)$ can be embedded in the polytope algebra
of  McMullen [12, 13, 14].  This is carried out by relating our
Minkowski weights, which  depend only on the lattice, to
McMullen's notion of weights on a polytope,
which depend on metric geometry in
Euclidean space.  In fact, we show  (Section 4) that,
the polytope algebra is the direct limit of all the Chow  rings,
taken with rational coefficients, as $ X$  varies over all
compactifications of a fixed torus. Our formula for multiplying
Minkowski weights is shown to be equivalent to a mixed volume
decomposition of McMullen (Proposition 4.3).

	Any variety has a Todd class in its Chow homology group with
rational coefficients.  If  $X$  is the toric variety of a simplicial fan,
this  Todd class is Poincar\'e  dual to a class $ Td_X$  in the Chow
cohomology, and hence to a  ${\bf Q}$-valued Minkowski
weight.  If $ v_1, . . . , v_d$
are the primitive lattice points along the edges of the fan, and
$ D_1,\ldots,  D_d $
 are the corresponding divisors, and  $ D = a_1 D_1 + . . . + a_d D_d $ is a
Cartier divisor whose line bundle is generated by its sections, the degree
of  $exp(D)\cdot Td_X$  is the number of lattice points $ u  \in  M  $
such that  $\langle u,v_i \rangle  \geq  -a_i  $ for
 $ 1 \leq i \leq d$.   When  $X$  is nonsingular, this is a
polynomial of degree $ n$  in the variables
$a_1, . . ., a_d$.  This is discussed  in Section 5.

We thank M.~Brion, V.~Batyrev,
R.~MacPherson, R.~Morelli, and B.~Totaro for
useful  conversations.  V.~Batyrev has pointed out that the spaces
$Hom(A_k(X),{\bf C})  $ for a complete toric variety
$X$ have appeared as the cohomology groups   $H^k(X,\Omega^k_X)$
(see [3, \S12] and [1]);  B.~Totaro's generalizations of
these results are discussed at the end of Section 1.

\vskip .3cm

\beginsection 1. Chow groups and the Chow cohomology ring

The Chow group $A_k (X)$ of an arbitrary variety $X$
is generated by all $k$-dimensional closed subvarieties
of $X$, with relations generated by divisors of rational
functions on $(k+1)$-dimensional subvarieties.
If $X$ is a toric variety, then there is a finite
presentation of $A_k(X)$ in terms
of torus-invariant subvarieties and torus-invariant
divisors. This is the content of Proposition 1.1 below.
Our notations concerning toric varieties are as in [6].

Let $X = X(\Delta)$ be the toric variety corresponding
to a fan $\Delta$ in a lattice $N$ of dimension $n$.
The torus-invariant closed subvarieties of $X$
are of the form $V(\sigma)$, as $\sigma$ varies
over the cones in $\Delta$, with $dim(V(\sigma))
= codim(\sigma) = n - dim(\sigma)$.
Each cone $\tau$ determines a sublattice $M(\tau)
= \tau^\perp \cap M $ of the lattice $M$ dual to $N$.
Each nonzero element $u \in M(\tau)$ determines a
rational function $x^u$ on $V(\tau)$. The divisor
of this rational function is
$$ [ div(x^u) ] \quad = \quad
\sum_\sigma \,\langle u , n_{\sigma,\tau} \rangle
\cdot [V(\sigma)] , \eqno (1.1) $$
where the sum is over all cones $\sigma$
that contain $\tau$ with $dim(\sigma) =
dim(\tau) + 1$, and $n_{\sigma,\tau}$ is a lattice
point in $\sigma$ whose image generates the
$1$-dimensional lattice $N_\sigma/N_\tau$.
Here $N_\sigma$ and $N_\tau$ are the sublattices
of $N$ generated by $\sigma \cap N$ and
$\tau \cap N$ respectively.

\proclaim Proposition 1.1.
\item{(a)} The Chow group $A_k(X)$ of a toric
variety $X = X(\Delta)$ is generated by the classes
$[V(\sigma)]$ where $\sigma$ runs over all cones of
codimension $k$ in the fan $\Delta$.
\item{(b)} The group of relations on these generators
is generated by all relations (1.1) where
$\tau$ runs over cones of codimension $k+1$ in
$\Delta$, and $u$ runs over a generating set of
$M(\tau)$.
\item{ }

\vskip -.4cm

\noindent {\sl Proof: }
Part (a) equals the proposition on page 96 in [6, \S 5.1].
Part (b) follows directly from Theorem 1 in [8]. \Box

\vskip .2cm

If $X$ and $Y$ are arbitrary varieties,
then there is a {\it K\"unneth map} $\,A_k(X)
\otimes A_l(Y) \rightarrow A_{k+l}(X \times Y)$, which is
defined by sending $[V] \otimes [W] $ to the class
of the product $[V \times W ]$. In the special
case where $X$ and $Y$ are toric, this map
is an isomorphism.

\proclaim Corollary 1.2.
If $X$ and $Y$ are toric varieties, then the
K\"unneth map $\, A_*(X) \otimes A_*(Y)
\rightarrow A_* ( X \times Y ) \,$ is an isomorphism.

\noindent {\sl Proof: \ }
This is a special case of Theorem 2 in [8]. In our toric
situation Corollary 1.2 may also be derived directly
from the presentation in Proposition 1.1, using the
fact that $\, V(\sigma \times \tau) \, = \,
 V(\sigma ) \times V(\tau) $.
\Box

\vskip .2cm

If $X$ is a non-singular toric variety, then the
Chow groups $A_k(X)$ are free abelian, and they are
never trivial for complete $X \, $ (see [6, \S 5.2]).
For general toric varieties these groups may
have torsion, and they may be trivial even if $X$
is complete.

\vskip .1cm

\noindent {\bf Example 1.3. } \
Let $\Delta$ be the fan over the faces of the cube
with vertices at $(\pm 1, \pm 1, \pm 1)$ in ${\bf Z}^3$.
Let $N$ be the lattice spanned by these vertices, that is,
$\, N \,= \,\{\,(x,y,z) \in {\bf Z}^3 \,:\, x \equiv y
\equiv z \,\,\,  (mod \, 2) \,\}$.  For any positive
integer $k$, let $\Delta_k $ be the
complete non-projective fan of the same combinatorial
type obtained from $\Delta$ by replacing the generator
$(1,1,1)$ by $(1,1,2k+1)$. An explicit computation using
Proposition 1.1 shows that $\,A_1 ( X( \Delta_k )) \, =
\,{\bf Z}/k {\bf Z}$. In particular, we have
$$ A_1(X(\Delta)) \, = \,{\bf Z} \,,\quad
A_1(X(\Delta_1 )) \, = \,0 \quad \hbox{and} \quad
A_1(X(\Delta_2 )) \, = \,{\bf Z}/2 {\bf Z}  .$$
\vskip -.4cm
\hfill \Box

\vskip .3cm

For some examples of Proposition 1.1 when $k=n-1$
see [6, \S 3.4]. If $\Delta$ is complete, then Proposition 1.1
implies that $\,A_0 (X(\Delta)) \, = \,{\bf Z}$.
The isomorphism from   $A_0(X(\Delta))$ to ${\bf Z}$ is
the degree map (to be considered
in (1.2) below). It satisfies $\,
deg\bigl( [V(\sigma)] \bigr)\,= \, 1 \,$ for all
$n$-dimensional cones $\sigma \in \Delta$.
If $\Delta$ is not complete, then
Proposition 1.1 implies that $\,A_0 (X(\Delta)) \, = \,0$.

When $X$ is a nonsingular variety, there is a natural
ring structure on the Chow homology $A_* (X) = \bigoplus_k A_k(X)$.
For a general singular variety $X$, there
is no natural ring structure on the  Chow homology, just as there is no
natural ring structure on the singular homology $H_*(X,{\bf Z})$.
Instead, there are {\it operational Chow cohomology groups}
$A^k(X)$. These groups, which were introduced by
the first author and R.~MacPherson in [7],
fit together to form a ring suitable for
doing intersection theory on $X$. An element
of $A^k(X)$ is a compatible collection of homomorphisms
of Chow groups
from $A_i (Y)$ to $A_{i-k}(Y)$ for all varieties
$Y$ mapping to $X$ and all integers $i \geq k$.
(See [5] or [7] for precise definitions and details.)
Composition of these homomorphisms defines the
multiplication $\, A^k(X) \otimes A^l(X) \rightarrow
A^{k+l}(X) \,$ which turns the Chow cohomology $A^*
(X) = \bigoplus_k A^k(X)$ into a graded commutative
${\bf Z}$-algebra.

For any complete variety $X$ there is a degree
homomorphism  $\,deg : A_0(X) \rightarrow {\bf Z}$.
An element of the Chow cohomology group $A^k(X)$
gives a homomorphism of Chow groups from
$A_k(X) $ to $A_0(X)$, and, by composition with ``$deg$'',
it gives a homomorphism from $A_k(X)$ to ${\bf Z}$. The
resulting  homomorphism of abelian groups is denoted
$$ {\cal D}_X \,: \, A^k(X) \quad \rightarrow
\quad Hom( A_k(X), {\bf Z}). \eqno (1.2) $$
For a general complete variety $X$ the map ${\cal D}_X$
can have a large kernel, even if $X$ is nonsingular.
On the other hand, it was shown in [8, Theorem 3] that
${\cal D}_X$ is an isomorphism if $X$ is a complete scheme
on which a connected solvable group acts
with finitely many orbits.
This class of schemes includes toric varieties and their
closed torus-invariant subschemes, whence
we get the following result.

\proclaim Proposition 1.4.
If $X$ is a complete toric variety, then the map
${\cal D}_X$ is an isomorphism. More generally,
if $Y$ is closed and torus-invariant in $X$,
then ${\cal D}_Y$ is an isomorphism.

 B.~Totaro has shown that, for a complex toric variety $X$,
 the canonical map from  $A_k(X)$ to Borel-Moore
 homology  $H^{BM}_{2k}(X)$  maps  $A_k(X)_{\bf Q}$
isomorphically onto the smallest weight space
$W_{2k}H^{BM}_{2k}(X)$.  He proves this in [17] for a class of
varieties containing all spherical varieties.  In the toric
case, Totaro gives the following simple description for these
weight spaces.  For any positive integer  $m$, let
$t_m : X \rightarrow X$  be the morphism determined by
multiplication by  $m$  on the underlying lattice  $N$.  Then
$W_{2k}H^{BM}_{2k}(X)$  is the subspace of  $H^{BM}_{2k}(X)$
of classes  $\alpha$  such that  $(t_m)_*(\alpha) = m^k\,\alpha$
for all $m$.  For  $X$  complete, this identifies $A^*(X)_ {\bf Q}$
with the corresponding weight space of  $H^*(X;{\bf Q})$, which
refines the result of Batyrev [1].

\vskip .3cm

\beginsection 2. Minkowski weights

Let $\Delta$ be a complete fan in a lattice $N$,
corresponding to a complete toric variety
$X = X(\Delta)$. Let $\Delta^{(k)}$
denote the subset of all cones of
codimension $k$.
An integer-valued function on $\Delta^{(k)}$
is called a {\it weight} of codimension $k$ on
$\Delta$. We say that $c$ is
a {\it Minkowski weight} if it satisfies the relation
$$ \sum_{ \sigma \in \Delta^{(k)} :\,
\sigma \supset \tau}
\!\!\! \langle u, n_{\sigma,\tau} \rangle \cdot
c(\sigma) \quad = \quad 0 \eqno (2.1) $$
for every cone $\tau $ in $\Delta^{(k+1)}$ and
every element $u$ of the lattice $M(\tau)$.
Equivalently, $c$ is a Minkowski weight if
$\sum_\sigma c(\sigma) n_{\sigma,\tau} $ lies
in $N_\tau$ for all $\tau \in \Delta^{(k+1)}$.
As in (1.1) the lattice point $n_{\sigma,\tau}$
is any representative in $\sigma$ for the
generator of the $1$-dimensional lattice
$N_\sigma/N_\tau$.  We have the following combinatorial
description of the Chow cohomology groups.

\proclaim Theorem 2.1.
The Chow cohomology group $A^k(X)$ of
a complete toric variety $X = X(\Delta)$
is canonically isomorphic to the group
of Minkowski weights of codimension $k$ on $\Delta$.

\noindent {\sl Proof: } This follows immediately
from combining Propositions 1.1 and 1.4. \Box

\vskip .1cm

\proclaim Remark 2.2. \rm Theorem 2.1 holds for every
closed and torus-invariant subscheme $Y$
of a complete toric variety $X$ as well.
The Chow cohomology of such a $Y$ is isomorphic to the
group of Minkowski weights that are supported on the cones
$\sigma \in \Delta $ with  $V(\sigma) \subset Y$.

\vskip .2cm

We shall rewrite the defining relation for Minkowski
weights of codimension $1$ on a complete fan $\Delta$.
Given cones $\rho \supset \sigma$ of codimension
$0$ and $1$ in $\Delta$, let $m_{\rho,\sigma}$ denote
the unique generator of the $1$-dimensional lattice
$M(\sigma)$ which is non-negative on $\rho$.
If $\tau $ is a codimension $2$  cone in $\Delta$,
then the star of $\tau$ in $\Delta$ is given by a sequence
of cones  $$ \sigma_1 \subset \rho_1 \supset \sigma_2
\subset \rho_2 \supset  \,\cdots\, \supset \sigma_t
\subset \rho_t \supset \sigma_{t+1} = \sigma_1 \eqno (2.2)
$$ with $codim(\sigma_i) = 1$
and $codim(\rho_i) = 0
$.

\proclaim Lemma 2.3.
A weight $c$ on the codimension $1$ cones in a complete fan $\Delta$ is a
Minkowski weight if and only if $\,\ \sum_{i=1}^t  c(\sigma_i)
\cdot m_{\rho_i,\sigma_i} \,=\,0 \,$ for all
cones $\tau$ of codimension $2$.

\noindent {\sl Proof: }
We identify both $M(\tau)$ and its dual lattice
$N / N_\tau$ with ${\bf Z}^2$. The primitive lattice vectors
$m_{\rho_i,\sigma_i}$ and $n_{\sigma_i,\tau}$ are orthogonal to
each other. There exists a consistent orientation around
the star of $\tau$, such that if
$m_{\rho_i,\sigma_i}$ has coordinates $(x_i,y_i)$, then
$n_{\sigma_i,\tau}$ has coordinates $(-y_i,x_i)$.
Being a Minkowski weight means that
$\sum_{i=1}^t c(\sigma_i) \cdot n_{\sigma_i,\tau} \,=\,0\,$
in $N/N_\tau$. This is equivalent to the condition
$\sum_{i=1}^t c(\sigma_i)\cdot m_{\rho_i,\sigma_i} \,=\,0\,$
in $M(\tau)$.\quad \Box

\vskip .1cm

For any variety $X$ there is a canonical homomorphism from
$Pic(X)$ to $A^1 (X)$. It takes a line bundle $L$
to the operator $\alpha \mapsto c_1(f^*L) \cap \alpha$
for $f : Y \rightarrow X$,
$\alpha \in A_* (Y)$. This can fail to be an isomorphism
[5, Ex.~17.4.9], but this does not happen for toric varieties.

\proclaim Corollary 2.4. {\rm (Brion, [2])} \ \  If $X$ is a
complete toric variety, then the canonical map  $Pic(X)
\rightarrow  A^1 (X)$ is an isomorphism.

This result was proved  for projective spherical varieties in
[2], and more recently Brion extended it to arbitrary complete
spherical  varieties. What follows is a simple proof for the toric
case.

\vskip .1cm

\noindent {\sl Proof: }
By Proposition 1.4 we may compose with the isomorphism
${\cal D}_X$ and show that the resulting map
$Pic(X) \rightarrow A^1 (X) \simeq
 Hom(A_1 (X), {\bf Z})$ is an isomorphism.
This map takes a line bundle $L$ to the linear functional
which assigns to a curve $C$ on $X$ the degree of $L$
on $C$. We know [6, \S 3.4] that the group
$Pic(X)$ of isomorphism classes of line bundles equals the
group of $T$-Cartier divisors modulo principal Cartier divisors.
To specify an element $L$ in $ Pic(X)$, we must define
elements $u(\rho)$ in $M$ for every $n$-dimensional cone
$\rho$ in the fan of $X$, such that whenever two cones
$\rho$ and $\rho'$ intersect in a codimension $1$ cone $\sigma$,
then $u(\rho) - u(\rho')$ lies in $ M( \sigma)$.
In this case there exists a unique integer $c(\sigma)$ such that
$$ u(\rho) - u(\rho') \quad = \quad
c(\sigma) \cdot m_{\rho,\sigma}  \quad = \quad
- c(\sigma) \cdot m_{\rho',\sigma} . \eqno (2.3) $$
The integer $c(\sigma)$ is the degree of $L$
on the invariant curve $C = V(\sigma)$.

We claim that the assignment $u \mapsto c$ represents
the map $Pic(X) \rightarrow A^1 (X)$.
Indeed, this assignment defines a
homomorphism of abelian groups from the
$T$-Cartier divisors to the weights of codimension $1$
on the fan of $X$.
The kernel of the map $u \mapsto c$
equals the group
$M$ of principal Cartier divisors, and
Lemma 2.3 verifies that its image lies in the
subgroup of Minkowski weights.

To complete the proof of Corollary 2.4, it suffices to
construct  the inverse map $A^1 (X) \rightarrow Pic(X) , \,
c \mapsto u$. Given any Minkowski weight $c$ of
codimension $1$, we  set $u(\rho_0) := 0$ for one fixed
$n$-dimensional  cone $\rho_0$. Then
we define $u(\rho)$ for all other cones $\sigma$
using the relation (2.3) along any path of
adjacent $n$-dimensional cones in the fan of $X$.
Since any two such paths with the same endpoints
differ by a sum of cycles like (2.2),
 Lemma 2.3 guarantees that this definition is
independent  of the chosen path. \Box

\vskip .2cm

\noindent{\bf Example 2.5. }\
{\sl (Minkowski weights on the hypersimplex)} \hfill \break
An important family of singular toric varieties arises
from the action of the torus $T = ({\bf C}^*)^n$
on the Grassmannian $ Gr_k ({\bf C}^n) $. We let
$ X_{k,n} $ denote the closure of the $T$-orbit of a
generic point in $ Gr_k ({\bf C}^n) $. This is an $(n-1)$-dimensional
projective toric variety, which is singular for
$ 2 \leq k \leq n-2 $. In what follows we describe its
Chow cohomology groups.

The polytope associated with
the Pl\"ucker embedding of $X_{k,n}$ is the {\it hypersimplex}
$$ \Delta(k,n) \quad = \quad conv \,
\{ \, e_{i_1} +  e_{i_2} +  \cdots + e_{i_k} \, : \,
1 \leq i_1 < i_2 < \cdots < i_k \leq n \,\} .$$
The complete fan of the toric variety
$X_{k,n}$ is the normal fan of $\Delta(k,n)$, considered with
respect to the lattice $ N = {\bf Z}^n / {\bf Z} (1,1,\ldots,1)$.
The hypersimplex $\Delta(k,n)$ is an important polytope
which appears in many different contexts;
for instance, see [9, \S 2.1].

We identify the cones of dimension $d$ in this fan
with the faces of codimension $d$ of $\Delta(k,n)$.
Here are a few basic facts about the hypersimplex:
\item{$ \bullet $} The hypersimplices
$\Delta(1,n)$ and $\Delta(n-1,n)$ are regular
$(n-1)$-simplices.
\item{$ \bullet $} The polytopes $\Delta(k,n)$
and $\Delta(n-k,n)$ are isomorphic.
Both have ${n \choose k} = {n \choose n-k}$ vertices.
\item{$ \bullet $}
For $ 2 \leq k \leq n-2 $, $\Delta(k,n)$ is an
$(n-1)$-dimensional polytope with $2n$ facets.
Its facet normals are the directed unit vectors
$\, \pm e_i \,$ considered modulo $span(e_1 + \ldots + e_n )$.

\vskip .1cm

\noindent
The positive-dimensional faces of $\Delta(k,n)$ are labeled by pairs $(I,J)$
where $ I , J \subset [n] := \{1,2,\ldots,n\}$,
$I \,\cap \, J = \emptyset $,
$|I| < k $, and $|J| < n - k $. The face with label $(I,J)$ equals
$$ {\cal F}_{I,J} \quad = \quad
conv \,\biggl\{ \, \, (\sum_{i \in I} e_i )\,
+ \, e_{\nu_1} + \cdots + e_{\nu_{k-|I|}} \,\,:\,\,
\{ \nu_1,\ldots,\nu_{k-|I|}\} \in { [n] \setminus
(I \cup J ) \choose
k - |I|} \,\biggr\}. $$
Thus the face ${\cal F}_{I,J}$ is affinely isomorphic
to the hypersimplex $\Delta( k - |I|, n - |I| - |J| )$.
In particular,  $\,codim({\cal F}_{I,J}) \, = \,
|I| +|J|$, and hence the  number of faces of
codimension $d$  of $\Delta(k,n)$ equals $$  f(d,k,n)
\quad = \quad  \sum_{i = max(0,k+d+1-n)}^{min(k-1,d)}
\pmatrix{ & n & \cr i \,,  & \! d-i \, , & \!\! n-d \cr}
\qquad \hbox{for} \quad 0 < d < n-1 . \eqno (2.4) $$

\proclaim Proposition 2.6.
The Chow cohomology group $ A^{ n - d-1}( X_{k,n} )$ is
isomorphic to the space of integer-valued functions
$c$ on the  faces of codimension $d$ of
$\Delta(k,n)$, which satisfy the following relations for
all codimension $d-1$ faces ${\cal F}_{I,J} \, $ and all
$\, r, s \in [n] \setminus (I \cup J ) $: \item{(a)}
$  c({\cal F}_{I \cup \{r\}, J}) +
 c({\cal F}_{I , J \cup \{s\}})
\,= \,
 c({\cal F}_{I \cup \{s\}, J}) +
 c({\cal F}_{I , J \cup \{r\}}) $,
if $|I| \! < \!  k -1 $ and $|J| < n-k-1 $;
\item{(b)}
$ c({\cal F}_{I , J \cup \{s\}})
\,\,= \,\,
c({\cal F}_{I , J \cup \{r\}}) $,
\ \ if $\,|I| = k -1 $ and $|J| < n- k-1 $;
\item{(c)}
$  c ({\cal F}_{I \cup \{r\}, J})
\,\,= \,\,
c ({\cal F}_{I \cup \{s\}, J}) $,
\ \ if $\,|I| <  k -1 $ and $|J| = n- k-1 $.

The proof of this proposition is straightforward using
the facets normals $\pm e_i \in N $ of $\Delta(k,n)$  and
the fact that all faces of
a hypersimplex are again hypersimplices.

\vskip .1cm

We do not know whether there exists a
nice general formula for the {\it Chow Betti numbers}
$\,\beta_{r,k,n} =  rank \,( A^{r}( X_{k,n} ))
=  rank \,( A_{r}( X_{k,n} ))$.
However, it is easy to see that
$\,\beta_{r,1,n} = \beta_{r,n-1,n} = 1 \,$ (because it is a
simplex), $ \beta_{0,k,n} = \beta_{1,k,n} = 1 \,$
(because all $2$-faces are triangles), and $\beta_{n-2,k,n} = n+1 $
for $2 \leq k \leq n-2$. Using direct computation
we also determined these numbers
\item{$\bullet$} for $ X_{2,4} $:
the Chow Betti numbers are $(1,1,5,1)$.
\item{$\bullet$} for $ X_{2,5}$: $\,\, (1,1,6,6,1) $.
\item{$\bullet$} for $ X_{2,6}$: $\,\, (1,1,7,7,7,1) $.
\item{$\bullet$} for $X_{3,6}$: $\,\, (1,1,7,22,7,1) $.

\noindent Generalizing the example of the hypersimplex, it would be
interesting to study the Chow cohomology of
the orbit closures of the torus action on a
general flag variety $G/P$ \ (see e.g.~[4]).
Do the Chow Betti numbers
of these toric varieties have combinatorial or
representation-theoretic significance~?
\Box

\vskip .3cm

Returning to the case of general toric varieties,
we shall briefly discuss the contravariant behavior of
the Chow cohomology in terms of Minkowski weights.
Let $\psi : N' \rightarrow N $ be a  homomorphism of
lattices, $\Delta$ a complete fan in $N$, and $\Delta'$
a complete fan in $N'$, such that each cone $\tau'$ in
$\Delta'$ is mapped under $\psi$ onto a subset of some
cone $\tau$ in $\Delta$. These data define
an equivariant morphism of complete toric varieties
$ f : X(\Delta') \rightarrow X(\Delta)$, and
conversely every such morphism arises in this way.
On the level of Chow cohomology there is an
induced ring homomorphism $ f^*$ from $A^*(X(\Delta))$ to
$A^*(X(\Delta'))\,$ [5, p.~324].

If $c$ is a Minkowski weight on $\Delta$, then $ f^* c$
is a Minkowski weight on $\Delta'$, and the problem is to
express $f^* c $ in terms of $c$.
In what follows we consider the special case where
$\psi \otimes {\bf Q} $ is surjective, or, equivalently,
where $f $ is dominant.
The general case, which is more difficult,
will be addressed in the next section.

\proclaim Proposition 2.7.
Let $ f  : X(\Delta') \rightarrow X(\Delta)$
be a dominant morphism of complete toric varieties
as above, and let $c \in A^k(X(\Delta))$ be a
Minkowski weight of codimension $k$ on $\Delta$.
Let $\tau' \in \Delta'$ with $codim(\tau') = k$,
and let $\tau $ be the smallest cone of $\Delta$
that contains $\psi(\tau')$. Then
$$ ( f^* c)(\tau') \quad = \quad
\cases { c( \tau) \cdot [ N :\psi(N') + N_\tau ] & if $
codim(\tau) = k $,\cr
          \quad 0    & if $ codim(\tau) < k $.  \cr} $$

\noindent {\sl Proof: }
The projection formula for Chow cohomology [5, p.~325] states that
$$ f_* ( f^* c\,\cap \,[V(\tau')])
\quad  = \quad  c \,\cap \, f_* ( [V(\tau') ] )
\qquad \hbox{in} \,\,\, A_0(X(\Delta)). \eqno (2.5) $$
By [5, Thm.~I.1.4] and [6, p.~56], we have
$$ f_*([V(\tau')])  \quad = \quad \cases{
[R(V(\tau')):R(V(\tau))]\cdot [V(\tau)]
& if $codim(\tau) = k $,\cr
  \quad 0    & if $ codim(\tau) < k $.  \cr} \eqno (2.6) $$
The degree of the field extension equals
$$ [R(V(\tau')):R(V(\tau))]
\,\,= \,\,[N/N_\tau : \psi_* (N' /  N_{\tau'})]
\,\,= \,\, [ N :\psi(N') + N_\tau ] , \eqno (2.7) $$
where $\psi_*$ is the induced homomorphism from
$N'/N_{\tau'}$ to $N/N_\tau$.
We substitute (2.7) into (2.6) and then into the
right hand side of (2.5). The assertion now follows
by applying the degree homomorphism to both sides of (2.5).
\Box

\vskip .3cm

\beginsection 3. Multiplication of Chow cohomology classes

In this section we describe
the ring structure of the  Chow cohomology
$A^*(X)$ of a complete toric variety $X=X(\Delta)$ in terms
of Minkowski weights on its fan $\Delta$.
We first recall the relevant results from [8].
Let $\gamma$ be a cone   in $\Delta^{(k)}\,$
and $V(\gamma)$ the corresponding
$k$-dimensional invariant subvariety of $X$.
We consider the diagonal embedding $\delta : X
\rightarrow X \times X $ and its restriction to
$V(\gamma)$. Using the isomorphism in Corollary 1.2,
we can write the diagonal in $V(\gamma)
\times V(\gamma)$ as a ${\bf Z}$-linear combination
$$ \bigl[ \delta ( \,V(\gamma) \,)\bigr] \quad
= \quad \sum_{\sigma,\tau} \,m^{\gamma}_{\sigma,\tau}\cdot
[ V_{\sigma} \times  V_\tau ]
\qquad \hbox{ in } \,\,\, A_k \bigl( \,V(\gamma) \times
V(\gamma) \,\bigr) ,  \eqno (3.1) $$
where the sum is over all pairs
$\sigma,\tau \in \Delta$ such that
$\gamma \subset \sigma$,
$\gamma \subset \tau$ and
$codim(\sigma) + codim(\tau) = codim(\gamma) = k $.
The coefficients  $m^{\gamma}_{{\sigma},{\tau}}$ are
generally not unique.
It was demonstrated in [8] that
knowledge of such coefficients characterizes both the
action of the Chow cohomology on the
Chow groups -- the {\it cap product} -- and
the multiplication within the Chow
cohomology -- the {\it cup product}.
In what follows we identify elements
of $A^*(X)$ with Minkowski weights on $\Delta$,
and we set $m^{\gamma}_{{\sigma},{\tau}} = 0$
if $\gamma \not\subset \sigma$ or
$\gamma \not\subset \tau$

\proclaim Proposition 3.1. \ {\rm ([8], Theorem 4) }
\item{(a)} If $c \in A^p(X)$ and
$\gamma \in \Delta^{(k)}$, then the
cap product $ c \cap [V(\gamma)] $
in $A_{k-p}(X)$ equals
$$  c \cap [V(\gamma)] \,\quad = \quad \!\!
\sum_{(\sigma,\tau) \in \Delta^{(p)}
\times \Delta^{(k-p)}} \!\!\!\!
m^{\gamma}_{{\sigma},{\tau}}\cdot c(\sigma) \cdot
[V(\tau)].$$
\item{(b)} If $c \in A^p(X)$ and $\tilde c \in A^q(X)$,
then their cup product $ c \cup {\tilde c} $ in $A^{p+q}(X)$
is the Minkowski weight given by the formula
$$ (c \cup {\tilde c}) (\gamma) \,\quad = \quad \!
\sum_{(\sigma,\tau) \in \Delta^{(p)}
\times \Delta^{(k-p)}} \!\!\!\!\!
m^{\gamma}_{{\sigma},{\tau}}\cdot c(\sigma) \cdot
{\tilde c}(\tau).$$

\vskip .1cm

Note that if $p=k$ in part (a)
then the sum on the right hand side reduces
to one term. Taking degrees on both sides we then get
$\,deg(c \cap [V(\gamma)]) \, = \, c(\gamma)$.

Proposition 3.1 shows that the intersection
theory on a toric variety $X = X(\Delta)$
is completely determined by a collection of
integers $\, m^{\gamma}_{{\sigma},{\tau}} \,$
satisfying formula (3.1). Our objective is
thus reduced to the problem of computing  $\,
m^{\gamma}_{{\sigma},{\tau}} \,$ for all $\sigma,
\tau,\gamma$ as above.
This problem is solved by the following theorem,
which we call the ``fan displacement rule''.

\proclaim Theorem 3.2.
If $\gamma$ is any cone in $\Delta$ and $ v $ a
generic element in the lattice $ N $, then
$$ [\,\delta(V(\gamma)) \,] \quad = \quad \sum_{\sigma,
\tau} m_{\sigma,\tau}^\gamma \cdot [V(\sigma) \times
V(\tau)],$$ where $m^\gamma_{\sigma,\tau} = [N :
N_\sigma + N_\tau]$, and the sum is over all pairs
$(\sigma,\tau)$ of cones in $\Delta$ such that  $\sigma$
meets $\tau + v $,  $\, codim(\sigma) + codim(\tau) =
codim(\gamma)$, and $\,\sigma,\tau \supset \gamma$.

Here ``generic'' means  outside a finite union of
proper linear subspaces in $N_\gamma \otimes {\bf R}$,
to be specified further below.
Note that the theorem stated in the introduction
is implied by Proposition 3.1 and Theorem 3.2.

\vskip .2cm

\noindent {\bf Example 3.3.}
We illustrate the fan displacement rule
for computing the cup product of Chow cohomology classes
in a simple example.
Let $X $ be the Hirzebruch surface $ {\bf F}_m$
 with fan $\Delta$
in $N = {\bf Z}^2$ having generators
$\,v_1 = (1,0),\,   v_2 = (0,1) , \,
   v_3 = (-1,m),\,  v_4 = (0,-1)$, where $m$ is a
nonnegative integer. Let $c$ and $\tilde c$ be
elements of $A^1(X)$, that is, Minkowski weights on the rays of
$\Delta$. The condition (2.1) around the zero-dimensional
cone $\tau = \{0\}$ states that
$$ c (v_1) = c(v_3) ,\,\,\,
c (v_2) + m \cdot c (v_3) = c(v_4) ,\,\,\,\,\,\,
\tilde c (v_1) = \tilde c(v_3) ,\,\,\,
\hbox{and} \,\,
\tilde c (v_2) + m \cdot {\tilde c} (v_3)= \tilde c(v_4).$$
We compute the cup product $\,c \cup \tilde c
\in A^2(X) = {\bf Z}^1 $, using
Theorem 3.2 and Proposition 3.1 (b).
It suffices to determine the number
$ \,(c \cup \tilde c) (\{0\})$.
We list six different choices
of a generic displacement vector $v \in N$,
giving six different but equivalent formulas:
\item{(a)} $\,\, v = (1,1)\,$
gives the formula $\,\,
(c \cup \tilde c) ( \{0\})
\,\, = \,\, c(v_1) \cdot {\tilde c}(v_4) \,+\,
            c(v_2) \cdot {\tilde c}(v_3) $.
\item{(b)} $\,\, v = (-1,1 +m )\,$
gives  $\,\,
(c \cup \tilde c) ( \{0\})
\,\, = \,\, c(v_2) \cdot {\tilde c}(v_1) \,+\,
            c(v_3) \cdot {\tilde c}(v_4) $.
\item{(c)} $\,\, v = (-1-m,+1)\,$
gives $\,\,
(c \cup \tilde c) ( \{0\})
\,\, = \,\, c(v_2) \cdot {\tilde c}(v_1) \,
         +\,c(v_3) \cdot {\tilde c}(v_2) \,
         +\,m \cdot c(v_3) \cdot {\tilde c}(v_1)  $.
\item{(d)} $\,\, v = (-1,-1)\,$
gives  $\,\,
(c \cup \tilde c) ( \{0\})
\,\, = \,\, c(v_3) \cdot {\tilde c}(v_2) \,+\,
            c(v_4) \cdot {\tilde c}(v_1) $.
\item{(e)} $\,\, v = (1,-1-m)\,$
gives  $\,\,
(c \cup \tilde c) ( \{0\})
\,\, = \,\, c(v_4) \cdot {\tilde c}(v_3) \,+\,
            c(v_1) \cdot {\tilde c}(v_2) $.
\item{(f)} $\, v = (1+m,-1)\,$
gives  $\,\,
(c \cup \tilde c) ( \{0\})
\, = \, c(v_1) \cdot {\tilde c}(v_2) \,
         +\,c(v_2) \cdot {\tilde c}(v_3) \,
         +\,m \cdot c(v_1) \cdot {\tilde c}(v_3)  $.
\Box

\vskip .2cm

In order to prove Theorem 3.2, we need a general lemma
for computing certain toric deformations on a toric
variety. Consider any toric variety $X = X(\Delta)$ of
dimension $n$, with $\Delta$ a fan in $N \simeq {\bf Z}^n$
and $T = T_N$ the dense torus in $X$. Let $L$ be a
saturated  $d$-dimensional sublattice of $N$. It
corresponds to a subtorus  $T_L$ of $T$. Let $Y$ be the
closure of $T_L$  in $X$. In what follows we determine the
class  $[Y]$ in the Chow homology group $A_d(X)$.
Our description will depend on the choice of an
element $v$ in $N$. For any $v \in N$ let
$$ \Delta(v) \quad := \quad
\{ \,\sigma \in \Delta \,:\,
L_{\bf R} + v \,\,\, \hbox{meets
$\,\sigma\,$ in exactly one point} \,\}. \eqno (3.2) $$
Note that $dim(\sigma) \leq n - d $ for all
$\sigma \in \Delta(v) $. We say that $v$ is
{\it generic} if $dim(\sigma) = n-d$ for all
$\sigma \in \Delta(v) $; in this case the unique
intersection point of $L_{\bf R} + v$ and $\sigma$
lies in the relative interior of $\sigma$.
Note that the set of generic
lattice points in $N$ is Zariski dense in $N_{\bf R}$.

\proclaim Lemma 3.4.
If $v $ is a generic lattice point in $ N$, then
$$  [Y] \,\, = \, \sum_{\sigma \in \Delta(v)} \!
m_\sigma \cdot [V(\sigma)] \quad \,\,\,
\hbox{in}  \,\,\, A_d(X), \qquad \hbox{ where }
\,\,m_\sigma \, = \, [N : L + N_\sigma] .$$

\noindent {\sl Proof of Lemma 3.4: \ }
We can assume that $\, L' = L + {\bf Z} v = L \oplus {\bf Z} v \,$
is a saturated sublattice of $N$. Indeed, since
$\Delta(v) = \Delta(v + \ell)$ and
$\Delta(r \cdot v) = \Delta(v)$ for all
$\ell \in L$ and $r \in {\bf N}$, we
can replace any given $v$ by an element of $N$
which is primitive modulo $L$.

We first consider the case $d=n-1$. By the discussion
in the previous paragraph, we may assume that
$N = L \oplus {\bf Z} v $. Let $M$ be the
lattice dual to $N$, and let $u$ be an element of
$M$ that is perpendicular to $L$ and satisfies
$\langle u,v \rangle = -1 $. Consider the rational
function $ \, f \, = \, x^u - 1 \,$ on $X$. To prove
Lemma 3.3 for $d=n-1$, it suffices to show
the following identity of Weil divisors on $X$:
$$  [ div(f) ] \quad = \quad [Y] \,\, - \,\,
\sum_{\sigma \in \Delta(v)}
m_\sigma \cdot [V(\sigma)] . \eqno (3.3) $$
Here $\,m_\sigma \,= \, [ N : L + N_\sigma ] \, = \,\,
min \, \{ m \in {\bf N} \,: \,
\exists \,\ell \in L \, \hbox{with} \,
mv + \ell \in \sigma \,\}$. The two sides of (3.3)
clearly have the same restriction $[T_L]$ to the torus $T$.
Hence it suffices to show that, for each ray  $\sigma$ of
$\Delta$, the irreducible $T$-divisor $[V(\sigma)]$
appears with the same coefficient on both sides of (3.3).
Note that $ord_{V(\sigma)}(f)$, the coefficient of
$[V(\sigma)]$ in $[div(f)]$, is a nonpositive
integer since $f$ has no zeros outside of $T$.

If $\sigma$ lies in $\Delta(v)$, then the unique point in
$\,(L_{\bf R} + v )\,\cap \,\sigma \,$ has the form
$\,v + {1 \over m_{\sigma}} \ell \,$ for some $\ell \in
L$. Then $  m_{\sigma} v + \ell \in N$ is the
primitive lattice point which generates $\sigma$, and
$$ ord_{ V(\sigma)} (f) \quad = \quad
 ord_{ V(\sigma)} (x^u) \quad = \quad
\langle  u , m_\sigma v + \ell \rangle \quad = \quad - m_\sigma ,$$
which is the required equation in this case.
If $\sigma$ does not lie in $\Delta(v)$, then the
generator of the ray $\sigma$ has the form $\ell - p v $ for
some $\ell \in L$ and some non-negative integer $p$.
If $p > 0 $ then $ ord_{ V(\sigma)} (x^u) = p > 0 $, so
$ord_{ V(\sigma)} (f) = 0 $, as required; if $p = 0$ then
$ x^u $ restricts to a nonconstant function on $V(\sigma)$,
so again $ ord_{ V(\sigma)} (f) = 0 $. This completes the
proof of Lemma 3.4 in the special case $d=n-1$.

For the general case we consider a saturated sublattice
$\, L' = L + {\bf Z} v = L \oplus {\bf Z} v \,$ of $N$.
Let $\Delta_{L'}$ denote the fan in $L'$ whose cones are the
intersections $\sigma \cap L'_{\bf R} \,$ for $\sigma
\in \Delta$. This gives a toric variety $X' =
X(\Delta_{L'})$ and a proper map $X'  \rightarrow X$.
Let $Y'$ be the closure of
the torus $T_{L}$ in $X'$. By the special case
proved above, we have an identity of divisor classes
$$ [Y'] \,\, = \,\, \sum_{\sigma '} m_{\sigma'}
[V( \sigma')] \qquad \hbox{in} \,\, A_d(X'), \eqno (3.4)
$$ where $\sigma'$ runs over all rays of $\Delta_{L'}$
that have the form $\sigma \cap L'_{\bf R}$
for some $\sigma \in \Delta(v)$.
The coefficients in (3.4) are
$$ m_{\sigma'} \quad = \quad
[ L'  : L + N_{\sigma'}] \quad = \quad
[ L' + N_\sigma : L + N_\sigma]. $$
To see that these two indices are equal,
apply the Second Isomorphism Theorem to the
quotient groups and use the fact that  $N_{\sigma'} =
N_{\sigma} \cap L'$.

Next we push the equation (3.4) forward from
$X'$ to $X$. The torus closure $Y'$ maps
birationally onto the torus closure $Y$, so $[Y']$ maps to
$[Y]$. Each $V(\sigma')$ maps to $V(\sigma)$, but the
degree of this mapping need not be one. To compute
its degree, we consider the restriction of this mapping
to the open torus $\,O_{\sigma'} $ in $ V(\sigma')$,
which maps onto the torus $\,O_{\sigma} $ in $
V(\sigma)$. The surjection of tori $\,O_{\sigma'}
\rightarrow O_{\sigma}\,$ is determined by the inclusion
of lattices $$ L' / N_{\sigma'} \,\, = \,\,
 L' / (L ' \cap N_{\sigma} ) \,\, = \,\,
(L' + N_{\sigma}) / N_\sigma \quad \hookrightarrow \quad
N / N_{\sigma} . $$
The degree of this map and hence of the map
$\,V(\sigma')  \rightarrow V(\sigma) \,$ equals
$[N : L'+ N_{\sigma}]$.

We conclude that the push-forward
of the equation (3.4) into $A_d(X)$ equals  $\,[Y] = $
$ \sum_{\sigma \in \Delta(v)} m_\sigma [V(\sigma)]$,
with coefficients
$$ m_{\sigma} \quad = \quad
 [ N : L'+ N_{\sigma} ] \cdot m_{\sigma'} \quad = \quad
 [ N : L'+ N_{\sigma} ] \cdot
 [ L'+N_{\sigma} : L+N_{\sigma} ]
\quad = \quad  [ N  : L+N_{\sigma} ]. $$
This completes the proof of Lemma 3.4. \Box

\vskip .2cm

We remark that, in the case when $X$ is projective,
Lemma 3.4 can also be derived from
Corollary 1.2 in [10].
The proof given there uses the method of Chow forms.
If $X$ is projective and $P$ the corresponding
polytope then $\Delta(v)$ is interpreted as the
{\it tight coherent face bundle} on $P$ defined by $L$ and $v$.

We are now prepared to prove our
fan displacement rule.

\vskip .2cm

\noindent {\sl Proof of Theorem 3.2.: }
We apply Lemma 3.4 to the diagonal embedding
$\, X \hookrightarrow X \times X$. This corresponds
to the diagonal inclusion of lattices
$\, \delta : N \hookrightarrow N \times N $.
Using a generic vector $(v,w) $ in $ N \times N $
gives the same formula in Lemma 3.4 as using the vector
$(v-w,0)$. Therefore it suffices to translate by
vectors $(v,0)$, for $v \in N$. Points of intersection
of the cone $\sigma \times \tau $ with
$\delta( N_{\bf R}) + (v,0)$ correspond to points
of intersection of $\sigma$ with $\tau + v $.
For a generic element $v \in N$ there is
at most one such intersection point whenever
$dim(\tau) + \dim(\sigma) = n$.
Lemma 3.4 implies that
$$ \delta_* ([X] ) \quad = \quad
\sum m_{\sigma,\tau} \cdot [V(\sigma) \times
V(\tau)], \eqno (3.5) $$
the sum over the pairs $(\sigma,\tau)$ such that
$\sigma $ meets $\tau + v$ in one point, with coefficient
$$ m_{\sigma,\tau} \quad = \quad
 [\, N \times N \,:\, \delta(N) + (N_\sigma
\times N_\tau) \,] \quad  = \quad
[\, N \, : \, N_\sigma + N_\tau \,]. \eqno (3.6) $$
The last equality is seen by mapping
$N \times N $ onto $N$ via $(a,b) \mapsto a- b$,
which has kernel $\delta(N)$ and maps
$N_\sigma \times N_\tau $ onto $N_\sigma + N_\tau$.
This proves Theorem 3.2 for the case $\gamma = \{0\}$.

In the general case we apply the formulas (3.5) and
(3.6) to the toric subvariety $V(\gamma)$ of $X$.
This means we work in the quotient lattice $N/N_\gamma$.
The multiplicity (3.6) remains the same: $$
m^{\gamma}_{\sigma,\tau} \quad = \quad [\, N/N_\gamma \, :
\, (N_\sigma + N_\tau)/ N_\gamma\,] \quad = \quad [\, N\,
: \, N_\sigma + N_\tau\,].$$ This completes the proof of
Theorem 3.2. \ \Box

\vskip .2cm

The Chow cohomology ring $A^*(X)$ operates
on the Chow groups $A_*(X')$ for every morphism
$f : X' \rightarrow X$. We shall describe this
operation for a toric morphism $f$.
Let $\psi : N' \rightarrow N $ be a  homomorphism of
lattices, $\Delta$ a complete fan in $N$, and $\Delta'$
a complete fan in $N'$, such that each cone $\sigma'$ in
$\Delta'$ is mapped under $\psi$ onto a subset of some
cone in $\Delta$.  If $\sigma$ is the smallest cone of $\Delta$
containing $\psi(\sigma')$, then the corresponding toric
morphism $\,f:  \,X'  =  X(\Delta') \,
\rightarrow \,X = X(\Delta) \,$ satisfies
$$ f \bigl( V(\sigma') \bigr) \quad \subseteq \quad
V(\sigma). \eqno (3.7) $$
Given an element $c \in A^k(X)$, i.e., a
Minkowski weight of codimension $k$ on $\Delta$, we wish to
describe the homomorphisms
$$ A_p( X') \,\, \rightarrow \,\, A_{p-k}( X' ),
\quad z \,\,\mapsto \,\, f^* c \,\cap \,  z. \eqno (3.8) $$
To this end we prove the following direct
generalization of Theorem 3.2.
Let $\,\delta_f : X' \rightarrow X' \times X
\,$ denote the graph of $f$.

\proclaim Theorem 3.5.
If $\gamma'$ is any cone in $\Delta'$ and
$v$ a generic element of $N$, then
$$ [ \delta_f( V(\gamma')) ] \quad = \quad
\sum_{\sigma',\tau} m^{\gamma'}_{\sigma',\tau}
[V(\sigma')] \otimes [V(\tau)]
\qquad \hbox{in} \,\,\, A_*(X' \times X) \,\, = \,\,
A_* ( X') \otimes A_* (X),$$
where the sum is over all pairs of cones
$\sigma' \in \Delta'$ and $\tau \in \Delta$
such that
$\gamma' \subset \sigma'$,
$\psi(\gamma') \subset \tau$,
and $codim(\sigma') + codim(\tau) = codim(\gamma')$, and
$$ m^{\gamma'}_{\sigma',\tau} \quad =\quad
\cases{  [N : \psi(N_{\sigma'}) + N_\tau]
& if $\,\psi(\sigma') \,$ meets $\tau + v$ \cr
\quad 0  & otherwise.\cr}$$

\noindent {\sl Proof: }
It suffices to prove this formula for
the special case $\gamma' = \{0\}$.
{}From this the general case is derived
by replacing $X'$ by $V(\gamma')$
and $X$ by $V(\gamma)$, where $\gamma$ is
the smallest cone of $\Delta$ containing $\psi(\gamma')$.
This is possible by equation (3.7).

We apply Lemma 3.4 to the
fan $\Delta' \times \Delta$ and the sublattice
$\, L = graph(\psi) \,$ of $N' \times N$,
with displacement vector $(0,-v) \in N' \times N$.
Clearly, the graph of $f$ is the closure in $X' \times X$
of the subtorus associated with $L$.
The translated lattice $\,L_{\bf R} - (0,v) \,$
meets a cone $\,\sigma' \times \tau\,$ in
$(N' \times N)_{\bf R}$ in exactly one point
if and only if $\psi(\sigma') $ meets $\tau + v$ in $N_{\bf R}$
in exactly one point. Since $v$ is generic,
this condition implies that
$\psi$ is injective on $\sigma'$ and
$dim(\psi(\sigma')) + dim(\tau) =
dim(\sigma') + dim(\tau) = n$.
Finally, the lattice index appearing in Lemma 3.4
simplifies to $[N : \psi(N_{\sigma'}) + N_\tau]$,
by an argument analogous to the derivation of (3.6).
\Box

\proclaim Corollary 3.6.
Let $c \in A^k(X)$. The homomorphism (3.8) is given by the
formula $$ f^* c \cap [V(\gamma')] \quad = \quad
\sum_{\sigma',\tau} \,m^{\gamma'}_{\sigma',\tau}\cdot
c(\tau) \cdot [V(\sigma')],$$
for any cone $\gamma' $ of codimension $p$ in $ \Delta'$.
Here the sum is over all pairs $\sigma',\tau$ as in
Theorem 3.5, subject to the additional condition that $\,
codim(\tau) = k \,$ and $\,codim(\sigma') = p-k$.

\noindent {\sl Proof: }
This follows immediately from Theorem 3.5 and
Corollary 1 in [8]. \Box

By considering the special case $p=k$ and
applying the degree homomorphism on both
sides, we obtain a general formula for
pullbacks of Minkowski weights.

\proclaim Corollary 3.7.
If $f : X' \rightarrow X$ is a morphism
of complete toric varieties, then the
homomorphism $\,f^* : A^k(X) \rightarrow
A^k(X')\,$ is given (in terms of Minkowski weights)
by the formula:
$$ (f^* c) (\gamma') \quad = \quad
\sum_{\sigma',\tau}\, m^{\gamma'}_{\sigma',\tau}\cdot
c(\tau) .$$
The sum is as in Corollary 3.6 subject to
$codim(\gamma') = codim(\tau) = k$, and
$codim(\sigma') = 0$.

In the special case were $\psi \otimes {\bf Q} $ is surjective (so
that $f$ is dominant) we recover Proposition 2.7 from Corollary 3.7.
In this case $\sigma'$ and $\tau$ are uniquely
determined by $\gamma'$ and $v$,
and the sum collapses to one term.

\vskip .3cm

\beginsection 4. Relations to the polytope algebra.

In this section we study the rational Chow
cohomology $\, A^*(X)_{\bf Q}
\,$ of a projective toric variety $X$. Our objective
is to relate $  A^*(X)_{\bf Q} $ to the polytope algebra
of McMullen [12, 13, 14]. Our main result in this section
(Theorem 4.2) expresses the polytope algebra as the direct
limit of the Chow cohomology rings of all compactifications
of a given torus, thus providing a cohomological
interpretation of McMullen's theory.

Structures similar to the polytope algebra
have been introduced also by Khovanskii-Pukhlikov [11]
and Morelli [15].
Here we restrict ourselves to McMullen's theory, with two minor
modifications: we work over ${\bf Q}$ instead of
${\bf R}$, and we replace $\Pi_0 \simeq {\bf Z}$
by ${\bf Q}$.

We first review the definitions.
The {\it polytope algebra} $\Pi$ is a ${\bf Q}$-algebra,
with a generator $[P]$ for every polytope $P$
in ${\bf Q}^n$, and $[\emptyset] = 0 $.
The generators satisfy the relations
\item{(V)} $\, [P \, \cup \, Q] \, + \,
[\, P \, \cap \, Q \,] \,\, = \,\, [P] \, + \, [Q] $,
whenever $P \,\cup \,Q $ is a polytope; and
\item{(T)} $\, [\,P + t \,] \,\,=\,\, [P] $, for all
translations $t \in {\bf Q}^n $.

\noindent
The multiplication in $\Pi$ is given by
\item{(M)} $\, [P] \cdot [Q] \,\, := \,\,
[P\,+ \,Q]$, where $P+Q = \{ p+q: p \in P, q \in Q \}$
is the {\it Minkowski sum}.

The multiplicative unit is the class of a point: $1 = [ \{0\} ]$.
A basic relation in $\Pi$ states that $([P] - 1)^{n+1} =
0$. This implies that the {\it logarithm}
of a polytope $P$ is well-defined:
$$ log([P]) \quad = \quad
\sum_{r=1}^n  {(-1)^{r+1} \over r} \,([P] - 1)^r . \eqno (4.1) $$
It is shown in [12] that $\Pi $ is a graded
${\bf Q}$-algebra, $\, \Pi \, = \,\bigoplus_{k = 0}^n \Pi_k$.
The $k$-th graded component $\Pi_k$ is the
${\bf Q}$-vector space spanned by all elements of the
form $(log([P]))^k $, where $P$ runs over all polytopes
in ${\bf Q}^n$.

We now fix a polytope $P \subset {\bf Q}^n$, and
we define $\Pi(P)$ to be subalgebra of $\Pi$
generated by all classes $[Q]$, where $Q$ is a
{\it Minkowski summand} of $P \,$
(i.e., $P \,= \, \lambda Q + R $, for some
positive rational $\lambda$ and some polytope $R$).
Let $\Delta$ denote the normal fan of $P$,
and let $X = X(\Delta)$ be the corresponding
projective toric variety. Here and throughout this section
we identify the lattice $N$ with the standard lattice
${\bf Z}^n$ inside ${\bf Q}^n$. Note that the algebra
$\Pi(P)$ depends only on the fan $\Delta$,
and hence it is an invariant of the toric variety $X$.

Every ample line bundle $D$ on $X$ gives rise
to a lattice polytope $P_D$ with normal fan $\Delta$.
More generally, if ${\cal O}(D)$ is generated by its
sections, then we get a polytope $P_D$ whose normal fan is
equal to or refined by $\Delta$. The latter condition is
equivalent to $[P_D] \in \Pi(P)$.
We have the following identification of the
polytope subalgebra  $\Pi(P)$ with a subalgebra
of the rational Chow cohomology of $X$.
The exponential of a divisor class $D$ on $X$ is
defined by the familiar formula: $\,
exp(D) = \sum_{r=0}^{dim(X)} {D^r / r !} $.

\proclaim Theorem 4.1.
There exists a monomorphism of graded ${\bf Q}$-algebras
$\, \theta \, : \, \Pi(P) \hookrightarrow A^*(X)_{\bf
Q}$ such that $\theta([P_D]) = exp(D)$ for every ample
divisor $D$ on $X$. The image of $\theta$ equals the
subalgebra of $ A^*(X)_{\bf Q}$ generated by the Picard
group $Pic(X) \otimes  {\bf Q} = A^1(X)_{\bf Q}  $.

If $P$ is a simple polytope,
or equivalently, if $\Delta$ is simplicial fan,
or equivalently, if $X$ is a $V$-manifold, then
the Picard group generates $A^*(X)_{\bf Q}$ as
an algebra. Theorem 4.1 implies that in this special case
the map $\theta$ is an isomorphism. This isomorphism was
already established by McMullen in [13, Thm~14.1].

\vskip .2cm

\noindent {\sl Proof: }
In the proof of Theorem 4.1 we shall
extend scalars and work over the
field of real numbers ${\bf R}$. The reason is that
some elementary analytic geometry arguments require
irrational constants. These irrationals drop out
in the final monomorphism $\theta \,$
(see also the discussion in [13, \S 15]).
We fix the standard inner product and Euclidean metric
on ${\bf R}^n$ and the induced inner product and metric on
each affine subspace of ${\bf R}^n$.

We recall the definition of weights given
in [13, \S 5]. Denote by ${\cal F}_k(P)$
the set of $k$-dimensional faces of the polytope $P$.
A {\it $k$-weight} on $P$ is a mapping
$\, \omega \,:\, {\cal F}_k (P) \rightarrow {\bf R} \,$
which satisfies the {\it Minkowski relations}:
$$ \sum_{F \subset G} \omega(F) \cdot v_{F,G}\quad
= \quad 0 \qquad \hbox{for all} \,\,\,
G \in {\cal F}_{k+1}(P), \eqno (4.2) $$
where the sum is over all $k$-faces $F$ of $G$,
and $v_{F,G}$ denotes the unit outer normal
vector (parallel to the affine span of $G$) to $G$ at its
facet $F$. The real vector space of all $k$-weights
on $P$ is denoted by $\Omega_k (P)$.

Every polytope $Q$ with $[Q] \in \Pi(P)$ defines
a $k$-weight $\omega$ as follows. We set
$\, \omega(F) \,:= \, vol_k (F') $, where
$F'$ is the face of $Q$ corresponding to the
face $F$ of $P$, and $vol_k ( \,\cdot \,)$ denotes
the standard $k$-dimensional volume form on
the affine span of $F'$.
Note that $\omega(F')$ may be zero if the normal
fan of $Q$ is strictly coarser than $\Delta$.
We write $\, \omega \, = \,
vol_k(Q) \,\in \,\Omega_k (P)$, and we call
$\omega$ the {\it $k$-th volume weight} of $Q$.
The fact that $\omega$ is indeed a weight
(i.e., that it satisfies (4.2)) is the content
of Minkowski's classical theorem.

It is shown in [13, Thm.~5.1] that the resulting map
$$ \phi \,:\, \Pi(P) = \bigoplus_{k=0}^n \Pi_k (P)
\,\rightarrow \, \Omega(P) := \bigoplus_{k=0}^n
\Omega_k(P), \qquad [Q] \,\mapsto \,\oplus_k \, vol_k(Q)
\eqno (4.3) $$ is a monomorphism of graded vector spaces.
The degree $k$ component of $[Q] =  exp \, (log\,[Q]) $
equals  $[Q]_k \, = \,{1 \over k\,!} (log\,[Q])^k $.
Therefore $\,\phi( log([Q])^k) \,= \,
k \, ! \cdot vol_k (Q)$.

We next show that $\Omega_k(P)$ is canonically
isomorphic to $ A^k ( X )_{\bf R}$,
the space of real-valued Minkowski weights on the
codimension $k$ cones in $\Delta$. Let
$\tilde M = M \otimes {\bf R}$,
$\tilde N = N \otimes {\bf R}$, and identify both spaces
with ${\bf R}^n$.
Suppose $F \in {\cal F}_k(P)$ and let $\sigma \in
\Delta^{(k)}$ be the normal cone to $P$ at $F$.
We identify $\tilde M (\sigma) = M(\sigma) \otimes {\bf R}$
with the  (vector space parallel to the) affine span
of $F$.  Let $\{ m_1, \ldots, m_k \}$ be an
orthonormal basis for $\tilde M (\sigma) $.
The orthogonal projection onto $\tilde M (\sigma)$
defines an isomorphism of vector spaces
$$ \eqalign{ proj_\sigma \quad : \quad &
\tilde N / \tilde N_\sigma \,\, \rightarrow \,\,
\tilde M (\sigma ) \cr
& v \, \mapsto \, \sum_{i=1}^k \,\langle v, m_i \rangle
\cdot m_i .\cr}  \eqno (4.4) $$

Let $Vol_\sigma$ denote the volume form
on $\tilde M (\sigma) $ which is normalized with
respect to the lattice $M(\sigma)$, i.e.,
$Vol_\sigma (T) = 1 $ for every primitive
lattice simplex $T $ in $M(\sigma)$. We
define a real constant $\nu_\sigma$ using
the ratio between the
normalized volume and the standard $k$-volume
$$\,  Vol_\sigma \,( \, \cdot \,) \quad = \quad
\nu_\sigma \cdot k \, ! \cdot vol_k \,( \, \cdot \,) .
\eqno (4.5) $$
Note that $\nu_\sigma$ is typically irrational and that
$0 < \nu_\sigma \leq 1 $.

We claim that the map
$$  \psi_k \,\,\, : \,\,\,
A^k (X)_{\bf R}  \,\,  \rightarrow \,\,
\Omega_k (P) ,  \quad \, \hbox{given by}
\quad \psi_k(c) (\sigma) \,:= \,
c(\sigma) / \nu_\sigma , \eqno (4.6) $$
is well-defined and is a vector space isomorphism.
To see this, suppose that $\tau \in \Delta^{(k+1)}$
is contained in $\sigma$.
The generator $n_{\sigma,\tau}$ of
$N_\sigma$ modulo $N_\tau$ satisfies
$$ || proj_\tau ( n_{\sigma,\tau} ) ||\quad = \quad
\nu_\tau  / \nu_\sigma . \eqno (4.7)$$
Let $c \in A^k(X)$ be any Minkowski weight.
Then the relation (2.1) translates into
$$  \sum_{\sigma \supset \tau}
{ \psi_k (c)(\sigma) } \cdot
{ proj_\tau (n_{\sigma, \tau}) \over || proj_\tau (
n_{\sigma,\tau}) || } \,\,\, = \,\,\,
{ 1 \over \nu_\tau} \cdot \sum_{\sigma \supset \tau}
c(\sigma) \cdot proj_\tau (n_{\sigma,\tau}) \,\,\, =
\,\,\, 0 \qquad \hbox{for all $\tau \in
\Delta^{(k+1)}$.}$$
This shows that $\psi_k (c)$
satisfies (4.2), with $G \supset F$
the faces of $P$ polar to $\tau \subset \sigma$.
Therefore $\psi_k (c)$ lies in $\Omega_k (P) $.
This argument also shows that $\psi_k$
is an isomorphism.

If we combine the maps $\psi_k$ for all $k$, then
we get a graded vector space isomorphism
$\, \psi \, :  A^* (X)_{\bf R}  \,
\simeq \, \Omega (P )$.
Let $\phi : \Pi(P) \hookrightarrow \Omega(P)$
be as in (4.3). The composition $\, \theta \,:= \,
\psi^{-1} \circ \phi \,$ is a
monomorphism of graded vector spaces
from  $\Pi(P)$ into $ A^* (X)_{\bf R}  $.

Let $D$ be any ample divisor on $X$, and $P_D$ the
corresponding lattice polytope in $\tilde M $.
The element $D^k$ of $A^k (X)$ is represented by the
Minkowski weight $\,\sigma \,\mapsto \,
Vol_\sigma (F) $, where $F$ is the $k$-face of $P_D$
polar to $\sigma \in \Delta^{(k)}$. This follows
from the corollary on page 112 of [6]. From (4.5) and
(4.6)  we derive $$ \psi( D^k) \quad = \quad
 k \, ! \cdot vol_k (P_D) \quad =  \quad \phi \bigl( \,
log ( [P_D])^k \bigr). \eqno (4.8) $$
This implies that $\, \theta ( log( [P_D] )^k ) \,= \,D^k
$, and therefore $\,\theta ([P_D]) \,= \,
exp(D)\,$ in $ A^* (X)_{\bf R} $.

The Picard group $ A^1 (X)_{\bf R}$ is spanned by
the ample divisors on $X$. This proves that the
degree $1$ component of the map $\theta$ is surjective,
and hence it is  a vector space isomorphism:
$$ \theta_1 \,\,: \,\, \Pi_1(P) \,\simeq \,
 A^1 (X)_{\bf R} \,,\,\,\, log([P_D]) \,\mapsto
\, D . \eqno (4.9) $$

The algebra $\Pi(P)$ is spanned as a
vector space by the classes $[P_D]$,
where $D$ runs over all ample divisors on $X$.
Our monomorphism $\theta$ is multiplicative
on these generators:
$$ \theta ( [P_D] \cdot [P_{D'}] ) \,\,\,
= \,\,\, \theta ( [P_D + P_{D'}] ) \,\,\,
= \,\,\, \theta( [P_{D + D'}]) \,\,\,
=\,\,\, exp ( D + D') \,\,= \,\, exp(D) \cdot exp(D').$$
In view of (4.9), this completes the proof of Theorem
4.1.
\Box

\vskip .2cm

Let  ${\cal Y}$ denote the family of all $n$-dimensional
complete toric varieties, that is, all compactifications
of an $n$-dimensional algebraic torus $T_N$. The
corresponding family of Chow cohomology rings
$\bigl\{ A^*(X) \bigr\}_{X \in {\cal Y}}$ is a
directed system: for every equivariant morphism
$f : X_1 \rightarrow X_2$ of $n$-dimensional
complete toric varieties there is an inclusion of
rings $\,f^* : A^*(X_2 ) \hookrightarrow
A^*(X_1) $; the fact that this is an inclusion follows from
the fact that the Chow cohomology groups are torsion free.
We can thus form the direct limit
$ \,\,\underline{\lim}_{\!{}_{>}}
A^*(X)  \,:= \, \,\,\underline{\lim}_{\!{}_{>}} \{
A^*(X) \}_{X \in {\cal Y}}$.
This limit of Chow cohomology rings
can be viewed as a  universal ring for
intersection theory on all compactifications of
a fixed torus. The following theorem shows that
this direct limit is an integral version of
the polytope algebra.

\proclaim Theorem 4.2.
The algebra $ \,\,\underline{\lim}_{\!{}_{>}} A^*(X)_{\bf Q}  \,$
is isomorphic to the polytope algebra $\Pi$.

\noindent {\sl Proof: } Consider a morphism
of $n$-dimensional  complete toric varieties $f : X_1
\rightarrow X_2$. Let $P_i$ be the
polytope associated with $X_i$, and
let $\theta^{(i)} : \Pi(P_i) \hookrightarrow
A^*(X_i)_{\bf Q}$ be the inclusion of Theorem 4.1.
Since $P_2$ is a Minkowski summand of $P_1$,
we have a natural inclusion $i : \Pi(P_2)
\hookrightarrow \Pi(P_1)$.
We claim that $\,f^* \circ \theta^{(2)}
\, = \,\theta^{(1)} \circ i $.
To see this, let $D$ be a divisor
on $X_2$ which is generated by its sections,
and let $P_D$ be the corresponding lattice
polytope. This defines a divisor $f^* D$ on
$X_1$ having the same polytope $\, P_{f^* D} = P_D $,
and therefore  $\,i([P_D]) \,= \,[P_D]\,=\, [P_{f^* D}]\,$
in $\Pi(P_1)$. In $A^*(X_1)_{\bf Q}$ we get the desired relation
$$  \theta^{(1)} (i([P_D]))
    \,\,= \,\,exp ( f^* D )
    \,\,= \,\, f^*  \bigl( exp ( D )  \bigr)
    \,\,= \,\, f^* \bigl( \theta^{(2)}([P_D])). $$
Here $f^*$  commutes with ``$exp$'' because it is
a ring homomorphism.

The polytope algebra $\Pi$ clearly equals the
direct limit of the finitely generated algebras $\Pi(P)$,
with respect to the inclusions $\Pi(P_2) \hookrightarrow
\Pi(P_1)$ whenever $P_2$ is a Minkowski summand of $P_1$.
The discussion in the previous paragraph
shows that there is a monomorphism $\theta$
from the polytope algebra $\Pi$ into
$ \,\,\underline{\lim}_{\!{}_{>}} A^*(X)_{\bf Q}  $.
It is compatible with all monomorphisms
$\theta^{(i)}$ as in Theorem 4.1. To see
that $\theta$ is surjective, it suffices
to consider equivariant resolution of
singularities: for every $X_2$ there exists
a morphism $f : X_1 \rightarrow X_2 $ such
that $X_1$ is smooth and hence $\theta^{(1)}$
is an isomorphism. This proves that
$\theta$ is an isomorphism.
\Box

\vskip .2cm

In Section 3 we saw that the cup product in $A^* (X)$
can be computed by a simple rule involving a generic displacement
of the fan $\Delta$. We shall now show that in
the projective case our rule is equivalent
to the following ``mixed volume computation'',
which was introduced by McMullen in [13, \S 5, p.~426].
Let $F$ and $G$ be polytopes in ${\bf R}^n$ such that
$dim(F+G) = dim(F) + dim(G) $. Then there exists a
unique real number $\alpha_{F,G}$, which depends
only on (the ``angle'' between) the affine
subspaces spanned by $F$ and $G$,
such that
$$ vol (F + G) \quad = \quad
\alpha_{F, G} \cdot vol (F) \cdot vol (G) ,
\eqno (4.10)$$
with respect to the standard volume forms in the
respective affine subspaces.
If $P$ is an $n$-dimensional polytope in $M_{\bf R} = {\bf R}^n$, then
each sufficiently generic linear functional $ v $ on $M$
defines a {\it mixed decomposition} $\Delta_v$
of the Minkowski sum $2P = P + P $. We recall
the definition of $\Delta_v$; for further details see
e.g.~[16]. Consider the polytope  $\,\widehat{2P} \, =\,
\{ (p + p', v(p)) \in {\bf R}^{n+1}  \, : \,
p,p' \in P \} $. A face of $\widehat{2P}$ is
a {\it lower face} it it has an
outer normal vector with negative last coordinate.
Each lower face projects onto a subpolytope
of $2P$ of the form $F+G$, where $F,G$ are
faces of $P$, and we have
$dim(F+G) = dim(F) + dim(G) $ by the genericity of
$ v $. The set of all such polytopes $F+G$
is a polyhedral decomposition of $2P$, which we
denote by $\Delta_v$ and call the {\it mixed
decomposition} of $P+P$ defined by $v $.

\proclaim Proposition 4.3. {\rm (McMullen, [14]) }
Let $x_1,x_2 $ be elements in the polytope algebra
$\Pi(P)$, given by their weights $\omega_i  = \phi(x_i) \in
\Omega(P)$. Then the weight of their product
$\omega = \phi(x_1 x_2) $ satisfies
$$  \omega(H) \quad = \quad \sum \alpha_{F,G} \cdot
\omega_1 (F) \cdot \omega_2(G) \, \qquad
\hbox{ for all } H \in {\cal F}_{p+q}(P) , \eqno (4.11) $$
where the sum is over all faces
$F,G$ of $H$ with $dim(F) + dim(G) = dim(H)$ and
$\,F+G \,\in \,\Delta_v$.

The fact that the formula (4.11) defines a weight,
i.e., that $\omega$ satisfies (4.2) whenever $\omega_1$
and $\omega_2$ do, is a non-trivial statement about polytopes.
This was proved by McMullen in [14]. It can also be  derived from
Theorem 4.1 and our results in Section 3, as follows:

\vskip .1cm

\noindent {\sl Proof:  }
Let $c \in A^p ( X)_{\bf R}, {\tilde c} \in A^q (X)_{\bf R}$ and
$c \cup {\tilde c} \in A^{p+q}(X)_{\bf R}$ their cup product, all
three represented by Minkowski weights on $\Delta$.
Our fan displacement rule in Theorem 3.2 states
$$ (c \cup {\tilde c}) (\gamma) \quad = \quad \sum
\,m^\gamma_{\sigma,\tau}  \cdot c (\sigma) \cdot {\tilde c} (\tau)
\, ,\qquad \hbox{for all $\gamma \in \Delta^{(p+q)}$},
\eqno (4.12) $$
where $\,m^\gamma_{\sigma,\tau} = [N : N_{\sigma} + N_{\tau} ]\,$
and the sum is over all pairs of
cones $\sigma \in \Delta^{(p)} $ and
$\tau \in \Delta^{(q)}$ such that
$\gamma \,\subset \, \sigma \,\cap \,\tau \,$
and $\,(\sigma + v) \,\cap \, \tau \,\not= \,
\emptyset $.

Let $F,G,H$ be the faces of $P$ corresponding to
$\sigma,\tau,\gamma$. Obviously, the condition $F,G \subset
H$ is equivalent to $\gamma \subset \sigma \cap \tau$.
The ``angle of linear
subspaces'' defined above satisfies
$$ \alpha_{F,G} \quad = \quad
m^\gamma_{\sigma,\tau} \,\cdot\, {\nu_\sigma \cdot
\nu_\tau  \over \nu_\gamma } . \eqno (4.13) $$
We set $ c = \psi(\omega_1)$,
$ {\tilde c} = \psi(\omega_2)$, where  $\psi$ is the isomorphism in (4.6).
Using this isomorphism and (4.13), we see that
(4.12) is equivalent to (4.11).
The proof of Proposition 4.3 is now completed by the
following lemma.

\vskip .1cm

\proclaim Lemma 4.4.
Let $F,G$ be faces of $P$ with normal cones
$\sigma,\tau $, and let $v \in  N$
be generic. Then $\,(\sigma + v) \,\cap \, \tau
\,\not= \, \emptyset \,$ if and only if $F+G$ is a face of
the mixed decomposition~$\Delta_v$.

\noindent {\sl Proof of Lemma 4.4: }
We have $F+G \in \Delta_v$ if and only if
$\,\widehat{F\!+\!G} \, = \,
\{ (p+q, v(p)) \in {\bf R}^{n+1} \,:\,
p \in F, q \in G \}\,$ is a lower face of $\widehat{2P}$
if and only if  there exists a linear functional
$\ell \in N_{\bf R} = ({\bf R}^n)^* $
such that $(\ell,1) \in ({\bf R}^{n+1})^* $
attains its minimum over $\widehat{2P}$
at $\widehat{F\!+\!G}$. This holds if and only if $\ell$
attains its minimum over $P$
at $G$, and $\ell+v$ attains its minimum
over $P$ at $F$, or, equivalently,
$\,\ell \in \sigma$ and $\ell +v \in  \tau $.
\Box

\vskip .3cm

\beginsection 5. The Todd weight of a smooth toric variety

The purpose of this section is to illustrate the role of Minkowski
weights in the context of a prominent application of toric varieties,
namely, counting lattice points (see [6, \S 5.3]
and the references given there).

Let $X$ be a nonsingular projective toric variety
of dimension $n$. In this situation the maps
$A_*(X) \rightarrow H_*(X)$ and
$A^*(X) \rightarrow H^*(X)$ are both isomorphisms,
if we take the ground field to be ${\bf C}$.
Therefore the isomorphism in  Proposition 1.4 reduces to the
Kronecker isomorphism $\,H^*(X) \simeq
Hom( H_*(X), {\bf Z}) \,$ from topology.
We recall from [5] that any variety $X$ has a Todd homology class
$td_X$ in $A_*(X)_{\bf Q}$. If $X$ is non-singular,
then $td_X = Td_X \cap [X]$, where $Td_X$ is the
Todd cohomology class  in  $A^*(X)_{\bf Q}$.

Let $\Delta $ be the unimodular fan of
$X$ in  $N \simeq {\bf Z}^n $.
By Theorem 2.1 the Todd cohomology
class $Td_X$ is presented by a
Minkowski weight $Td_\Delta$ on $\Delta$. We call
$ Td_\Delta $ the {\it Todd weight}. It is our objective to
express $ Td_\Delta $ in terms  of a
certain  multivariate Ehrhart polynomial $\Phi$.

Let $ v_1,\ldots,v_d  $ denote the primitive lattice
points in $N$ along the rays of $\Delta$.
Each $v_i$ corresponds to a divisor $D_i$ on $X$.
Let ${\cal K} (\Delta)$ be the family of all
lattice polytopes $P$ such that $\Delta$ equals
or refines the normal fan of $P$.
The polytopes in ${\cal K}(\Delta)$ are of the form
$$ P \quad = \quad
\{  {\bf x} \in {\bf R}^n \, | \,
{\bf x} \cdot { v}_1 \geq - a_1, \ldots ,
{\bf x} \cdot { v}_d \geq - a_d \,\}, \eqno (5.1) $$
where $(a_1,\ldots,a_d)$ runs over all lattice points
in a certain closed convex cone in ${\bf R}^d$.
They are in bijection with
the divisors $\,D = a_1 D_1 + \ldots + a_d D_d \,$
on $X$ whose line bundles are generated by their sections.
Let $\Phi = \Phi(a_1,\ldots,a_d)$ denote the number
of lattice points in the polytope $P$ in (5.1).
We are interested in this number as a function
of the parameters $a_1,\ldots,a_d$.

\proclaim Proposition 5.1.
\item{(a)} The function $\Phi $ is
polynomial of degree $n$ in the variables
$a_1,\ldots,a_d$. \item{(b)} If a monomial $ a_1^{e_1}
a_2^{e_2}\cdots a_d^{e_d} $ appears in the expansion of
$\Phi$, then the generators $v_j$ that occur with
positive exponent $e_j > 0 $ span a cone in $\Delta$.

\noindent
{\sl Proof: }
The rational Chow ring of $X$ equals
$ A^*(X)_{\bf Q} = {\bf Q}[x_1,\ldots,x_d]/I $,
where $I$ is the ideal generated by all
linear relations
$\,\sum_{i=1}^d \langle m ,  v_i \rangle
\cdot x_i \,$ where $  m  \in  M $, and all
square-free monomials $ x_{i_1} x_{i_2} \cdots x_{i_r}\,$
such that $ v_{i_1},\ldots,v_{i_r} $ do
not span a cone in $\Delta$.
Let $Td_X$ denote the Todd cohomology class of $X$.
By the Hirzebruch-Riemann-Roch Theorem, we have
$$ \Phi \quad =
\quad \int exp\,(D)\cdot Td_X  \quad = \quad
\sum_{i=0}^n {1 \over i \,!} \int D^i \cdot Td_X . \eqno
(5.2) $$ This formula implies part (a) of Proposition~5.1
because $\int D^i \cdot Td_X $ is a homogeneous polynomial
of degree $i$ in $a_1,\ldots,a_d$.
To prove part (b) we consider the expansion
$$ D^i \quad = \quad (a_1 x_1 + \ldots + a_d x_d )^i
\quad = \quad \sum_{\nu_1,\ldots,\nu_i} a_{\nu_1} \ldots
 a_{\nu_i} x_{\nu_1} \ldots x_{\nu_i}.$$
After reducing modulo the square-free monomial relations in $I$,
the sum on the right hand side contains only monomials
$ a_{\nu_1} \cdots a_{\nu_i}$ which are supported on
$\Delta$.
The Todd class $Td_X$ is a polynomial in
$x_1,\ldots,x_d$ with rational coefficients.
Therefore the expression (5.2) is
a ${\bf Q}$-linear combination of monomials
$ a_{\nu_1} \ldots a_{\nu_i}$ which are supported
on $\Delta$.
\ \ \Box

\vskip .2cm

We remark that Proposition 5.1 can also be derived from the
results in [11]. For a generalization of the polynomial
$\Phi$ to arbitrary complete schemes
see [5, Example 18.3.6]. In what follows we identify $ A^*
(X)_{\bf Q}$ with the ring of ${\bf Q}$-valued Minkowski
weights on $\Delta$, and as well as with the
quotient of $ {\bf Q}[x_1,\ldots,x_d]$ modulo $I$.
The following easy lemma explicates the isomorphism
between the two graded algebras.

\vskip .1cm

\noindent {\bf Lemma 5.2. } {\sl
Let $p$ be a homogeneous element of degree $i$
in ${\bf Q}[x_1,\ldots,x_d] $ representing
a class in $A^i(X)_{\bf Q}$. Then the
corresponding Minkowski weight equals
$$ c \,: \, \Delta^{(i) } \rightarrow {\bf Q} \, , \,\,
\quad \sigma \,\mapsto \,
\int \! p(x_1,\ldots,x_d) \cdot \!\! \prod_{ v_j \in
\sigma} x_j  . \eqno  (5.3) $$
}

\vskip .2cm

We are now prepared to relate
the polynomial $\Phi$ to
the Todd weight $Td_\Delta$.

\proclaim Theorem 5.3.
If $\sigma$ is a cone of codimension $i$
in the fan $\Delta$, then $\,Td_\Delta(\sigma)\,$
equals the coefficient of the square-free monomial
$\,\prod_{ v_j \in \sigma} a_j \, $  in the
polynomial $\Phi(a_1,\ldots,a_d)$.

\noindent
{\sl Proof: }
Let $Td^i_X$ denote the
degree $i$ component of the Todd class,
represented by a polynomial of degree
$i$ in $ x_1,\ldots,x_d$.
The coefficient of interest can be
computed by applying the differential operator $\,\prod_{{\bf
v}_j \in \sigma} \! { \partial \over \partial a_j }\,$ to the
degree $n-i$ component of the polynomial
$\Phi(a_1,\ldots,a_d)$. By (5.2), the degree $n-i$ component
equals $\, \int {1 \over (n-i)~!}
 (a_1 x_1 + \ldots + a_d x_d)^{n-i} Td^{i}_X $.
By differentiating under the integral sign, we get precisely
the expression in (5.3) for $p = Td^i_X$.
This proves the claim. \ \ \Box

\vskip .1cm

\proclaim Corollary 5.4.
The lattice point enumerator
$\Phi(a_1,\ldots,a_d)$ is uniquely determined
by the coefficients of its square-free terms.

\noindent {\sl Proof: }
This follows from Theorem 5.3 and the formula (5.2)
because the Todd class $Td_X$
is uniquely determined by the
values of the Todd weight $Td_\Delta$. \Box

\vskip .2cm

Consider the polytope $P$ as in (5.1) and the corresponding
divisor $D$. The Minkowski weight corresponding
to $exp(D)$ is denoted $[P]$ and called the
{\it volume weight}. This notation is consistent
with McMullen's polytope algebra.  By our results in the
previous section, the volume weight satisfies
 $\, [P](\sigma) \,= \, Vol_\sigma(P^\sigma)\,$
for all $\sigma \in \Delta$.
Thus the formula (5.2) can be rewritten as
$\,\Phi \,= \,  \int [P] \cdot Td_\Delta $.
This means that the polynomial $\Phi(a_1,\ldots,a_d)$
can be computed by the  ``fan displacement
rule'' presented in Section 3, namely, by multiplying
the Todd weight $Td_\Delta$
with the volume weight $[P]$. We illustrate
this in a small example.

\vskip .2cm

\noindent {\bf Example 5.5.} \ ($n=2,d=4$, {\sl Hirzebruch surfaces}). \  \
Consider the two-dimensional fan $\Delta$ with four
generators $v_1 = (1,0)$,  $v_2 = (0,1)$, $v_3 = (-1,m)$
for a nonnegative integer  $m$, and $v_4 =  (0,-1)$.
Then $P$ is the (possibly degenerate) quadrangle defined by the
inequalities  $\, x \geq - a_1, \,a_4 \geq y \geq - a_2 $, and
$\,-x + m y \geq - a_3$. To ensure that $P$ has the
correct normal fan (equivalently, $P \in {\cal K}(\Delta)$),
we need to assume $\, a_1 + a_3 \geq  m a_2$ and $ a_2 + a_4 \geq 0 $.

We wish to count the number $\Phi(a_1,a_2,a_3,a_4)$
of lattice points in $P$.
The volume weight $[P]$ equals
$$  \eqalign{
& \{0\} \,\,\, \mapsto \,\,\, a_1 a_2 + a_2 a_3 + a_3 a_4
+ a_4 a_1 - {m \over 2} a_2^2 + {m \over 2} a_4^2, \cr
& \{ { v}_1  \} \, \mapsto \, a_2 + a_4 ,\,\,
\{ { v}_2  \} \, \mapsto \, a_1 + a_3 - m a_2 ,\,\,
\{ { v}_3  \} \, \mapsto \, a_2 + a_4 ,\,\,
\{ { v}_4  \} \, \mapsto \, a_1 + a_3 + m a_4 ,\cr
&  \{ { v}_1,{ v}_2 \} \, \mapsto \,1 , \,\,
\{ { v}_2,{ v}_3 \} \, \mapsto \,1 , \,\,
\{ {v}_3,{ v}_4 \} \, \mapsto \,1 , \,\,
\{ { v}_1,{ v}_4 \} \, \mapsto \,1 . \cr} $$
The Todd weight $Td_\Delta $ equals
$$  \eqalign{ & \{0\} \, \mapsto \,1 \, , \,\,\,
 \{ { v}_1  \} \, \mapsto \, 1 , \,\,
\{ { v}_2  \} \, \mapsto \, 1 - m / 2 , \,\,
\{ { v}_3  \} \, \mapsto \, 1 , \, \,
\{ { v}_4  \} \, \mapsto \, 1 + m / 2 \cr
&  \{ {v}_1,{ v}_2 \} \, \mapsto \,1 , \,\,
\{ { v}_2,{ v}_3 \} \, \mapsto \,1 , \,\,
\{ {v}_3,{v}_4 \} \, \mapsto \,1 , \,\,
\{ { v}_1,{ v}_4 \} \, \mapsto \,1 . \cr} $$
Applying the ``fan displacement rule'' to
these two weights on $\Delta$, we get the formula
$$ \Phi \quad = \quad
 a_1 a_2
+ a_2 a_3
+ a_3 a_4
+a_1 a_4
-{m \over 2} a_2^2
+{m \over 2} a_4^2
+a_1
+ ( 1 - { m \over 2}) a_2
+ a_3
+ (1 + {m \over 2}) a_4
+ 1.$$
\hfill \Box

\vskip .1cm

In closing we comment on the general case when $X$ is singular.
The preceding discussion extends to the case where $\Delta$
is an arbitrary complete simplicial fan. On still has a
Todd cohomology class $Td_X$ in $A^*(X)_{\bf Q}$, and the
degree of $exp(D) \cdot Td_X$ counts the number of
lattice points in the polytope $P$, but only when
$D = \sum a_i D_i$ is a {\it Cartier divisor}.
For singular varieties, including this simplicial case,
the lattice point enumerator $\Phi$ is generally not
a polynomial but only a quasi-polynomial.
If the fan is not simplicial, moreover, there may be
no Todd weight at all, as in the following example.

\vskip .2cm

\noindent {\bf Example 5.6.} \
{\sl A three-dimensional complete projective
toric variety $X = X(\Delta)$ whose Todd homology class
$td_X $ is not in $A^*(X)_{\bf Q} \cap [X]$.
Thus there is no Todd weight on $\Delta$.}

\vskip .1cm
\noindent
The Todd homology class of a toric variety
satisfies $td_{n-1}(X) = {1 \over 2} \sum_{i=1}^d [D_i]$,
where the $D_i$ are the divisors corresponding to the rays
of the fan; see [6, \S 5.3]. Thus it suffices to find
a projective toric variety $X$ such that
$\sum_i [D_i]$ is not in the image of
$A^1(X)_{\bf Q} = Pic(X)_{\bf Q}$ under taking the cap
product with the fundamental class $ [X]$.
This is equivalent to saying that $\sum_i [D_i]$
is not a ${\bf Q}$-Cartier divisor, or that there
is no piecewise linear function $\psi$ on $N$ whose
values on all the primitive generators $v_i$ of the rays
of the fan are equal.

Let $P$ be the pyramid in $M = {\bf Z}^3$ with
apex $(0,0,1)$ and a random quadrilateral base,
say, with vertices $(2,1,-1)$, $ (1,-1,-1 )$,
$ (-3,-2,-1 )$, and $ (-1,1,-1 )$.
Let $\Delta$ be the normal fan of $P$ and
$\sigma \in \Delta $ the normal cone at the apex.
The primitive generators of $\sigma$ are
$(-4,2,-3)$,$(0,-2,-1)$, $(6,-4,-5)$, and $(-2,8,-5)$
in $N = {\bf Z}^3$. There is no nonzero element
$u $ in $M$ with equal values on these four vectors,
since
$$ det \, \pmatrix{
 -4 & \phantom{-} 2 & - 3  & 1 \cr
\phantom{-} 0 &  -2 & -1  & 1\cr
\phantom{-} 6 &  -4 & -5 & 1\cr
  -2  & \phantom{-} 8 & -5 & 1 \cr} \quad
= \quad 176 \quad \not= \quad 0 .$$

\vskip 1cm

{\baselineskip=12pt

\centerline{\bf References}

\vskip .2cm

\item{[1]} V.~Batyrev: ``Rational equivalence
and K-theory of toric varieties'' (in Russian), in
{\sl Algebra, Logic and Number Theory},
Moscow State University, 1986, no.~3, pp.~20--23.
\vskip .1cm

\item{[2]} M.~Brion: ``Groupe de Picard et nombres
charat\'eristiques des variet\'es sph\'eriques'',
{\sl Duke Math.~J.} {\bf 58} (1989) 397--425.
\vskip .1cm

\item{[3]}  V.I.~Danilov: ``The geometry of toric varieties'',
{\sl Russian Math.~Surveys} 33:2 (1978), 97--154.

\item{[4]} H.~Flaschka and L.~Haine:
``Torus orbits in G/P'', {\sl Pacific J.~Math.}
{\bf 149} (1991) 251--292.
\vskip .1cm

\item{[5]} W.~Fulton: {\sl Intersection Theory},
Springer Verlag, 1984.
\vskip .1cm

\item{[6]} W.~Fulton: {\sl Introduction to Toric Varieties},
Princeton University Press, 1993.
\vskip .1cm

\item{[7]} W.~Fulton and R.~MacPherson:
``Categorical Framework for the Study of Singular Spaces'',
Memoirs Amer.~Math.~Soc. {\bf 243}, 1981.
\vskip .1cm

\item{[8]} W.~Fulton, R.~MacPherson, F.~Sottile and
B.~Sturmfels:  ``Intersection theory on spherical varieties'',
{\sl Journal of Algebraic Geometry}, to appear.
\vskip .1cm

\item{[9]} I.M.~Gel'fand and R.D.~MacPherson:
``Geometry in Grassmannians and a generalization of the
dilogarithm'', {\sl Advances in Math.} {\bf 44} (1982) 279--312.
\vskip .1cm

\item{[10]} M.~Kapranov, B.~Sturmfels and A.~Zelevinsky:
``Quotients of toric varieties'',
{\sl Mathematische Annalen} {\bf 290} (1991), 643-655.
\vskip .1cm

\item{[11]} A.G.~Khovanskii and A.V.~Pukhlikov: ``Finitely
additive measures of virtual polytopes'', {\sl St.~Petersburg
Math.~J.} {\bf 4} (1993) 337--356.
\vskip .1cm

\item{[12]} P.~McMullen: ``The polytope algebra'',
{\sl Advances in Math.} {\bf 78} (1989) 76--130.
\vskip .1cm

\item{[13]} P.~McMullen: ``On simple polytopes'',
{\sl Inventiones math.} {\bf 113} (1993) 419--444.
\vskip .1cm

\item{[14]} P.~McMullen: ``Separation in the polytope algebra'',
{\sl Beitr\"age zur Geometrie und Algebra} {\bf 34} (1993) 15--30.
\vskip .1cm

\item{[15]} R.~Morelli: ``A theory of polyhedra'',
{\sl Advances in Math.} {\bf 97} (1993) 1--73.
\vskip .1cm

\item{[16]} B.~Sturmfels: ``On the Newton polytope of the
resultant'', {\sl J.~Algebraic~Combinatorics} {\bf 3} (1994)
207--236.

\item{[17]} B.~Totaro: ``Chow groups of linear varieties'',  preprint.

}
\bye